# How to Apply 3D³ Prediction? A Novel Mathematical Model to Generate Pareto Optimal Clinical Applicable IMRT Treatment Plan On the Foundation of Dose Prediction and Prescription

**Ali Yousefi[1], Saeedeh Ketabi[1*], Iraj Abedi[2]**

[1] Department of Management- Operations Research, Faculty of Economics and Administrative Sciences, University of Isfahan, Isfahan, Iran

[2] Department of Medical Physics, School of Medicine, Isfahan University of Medical Sciences, Isfahan, Iran

[*] Corresponding Author

**Abstract**

In this paper knowledge based planning has been revolutionized via a novel mathematical model which converts three dimensional dose distribution (3D³) prediction to a clinical utilizable IMRT treatment plan. Presented model has benefited from both prescribed dose as well as predicted dose and its objective function includes both quadratic and linear phrases, so it was called *QuadLin model*. The model has been run on the data of 30 patients with head and neck cancer randomly selected from the Open KBP dataset. For each patient, there are 19 sets of dose prediction data in this database. Therefore, a total of 570 problems have been solved in the CVX framework and the results have been evaluated by two plan quality approaches: 1- DVH points differences, and 2- satisfied clinical criteria. The results of the current study indicate a strong significant improvement in clinical indicators compared to the reference plan of the dataset, 3D³ predictions, as well as the results of previous researches. Accordingly, on average for 570 problems and total ROIs, clinical indicators have improved by more than 21% and 15% compared to the predicted dose and previous research, respectively.

**Keywords:** Three dimensional dose distribution (3D³) prediction; Head and neck cancer; Knowledge-based planning; Open KBP dataset; CVX framework

## I. Introduction

Cancer radiation therapy planning in IMRT technique includes two main issues: 1- Selecting the number of beams and optimal angles 2- Optimizing the dose delivered to the patient's voxels. The second issue, which is more important, involves the optimization of beam weights to define a fluence map for each radiation field. In other words, the intensity of each of the beamlets that make up the beam is optimized so that the three-dimensional dose distribution of the patient's voxels maximizes compliance with the prescribed dose as well as the standard thresholds of organs at risk. In the last decade, extensive researches have been conducted to facilitate the trial and error process in the radiation therapy planning process and techniques have been developed to eliminate IMRT programming defects, which are classified into the following three categories: 1- Knowledge Based Planning (KBP), 2- Multi-criteria optimization (MCO), and 3- Hierarchical optimization with limitations [1]. The main idea in Knowledge-Based Planning is to extract a model for predicting appropriate planning goals from a previously planned patient database. KBP uses a database of previously treated patients to create a model that matches the dose delivered to the patient's voxels with the final dose-volume histogram. Zarepisheh et al. (2014) proposed an IMRT optimization algorithm based on dose-volume histogram (DVH) for automated treatment planning and adaptive radiotherapy scheduling [2]. Momin et al. (2021) provided a survey on KBP methods during the last decade which classified into two major categories consisted of traditional KBP methods and deep-learning methods (see that paper for more review on KBP) [3]. In recent years, with the development of artificial intelligence and deep learning techniques, it has become possible to predict the three-dimensional distribution dose (3D³) of a new patient based on the treatment plans of similar recent patients. Therefore, some new questions arise for the above issue: how to make use of the predicted 3D³ obtained from deep learning to facilitate treatment planning? How to convert predicted 3D³ into clinical usable Pareto optimal plan? Little

---

researches have been done in this regard. Fan et al. (2019) attempted to produce treatment plan closer to the predicted dose by applying a simple norm 2 objective function trying to pull dose in each voxel to its predicted value [4]. Their study did not pay attention to how well the dose prediction was done. As a result, in the case of poor quality dose prediction, the created treatment plan is easily affected and becomes ineffective. Babier et al. (2020) entered the prescribed dose into the model very poorly, but still insist on achieving the predicted dose [5]; Therefore, in their research, the quality of the predicted dose was not given much importance. Zhong et al. (2021) developed a 3D-Uet deep learning model and then generated rectum IMRT plans using post-optimization strategies, including clinical dose target metrics homogeneity index (HI) and conformation index (CI) in order to make the generated plan closer to the reference plan. However, their results showed that manual plans perform slightly better than plans with post-optimization strategies, but they claim these plans are still within clinical requirements [6]. Babier et al. (2021) did a valuable job of presenting a database including data from 100 head and neck cancer patients, along with the dose predicted by the 21 participated teams in the Open KBP Grand Challenge [7]. Babir et al. (2022) presented four optimization models with dose mimicking approach and strengthened the role of the prescribed dose. The results of the recent study showed an increase of about 6% in the satisfaction of clinical indicators compared to the predicted dose [8].

In the current paper a novel mathematical model has been developed in order to generate Pareto optimal IMRT treatment plan using a complex of 3D³ prediction and dose prescription values. Some dosimetry factors and constraints such as maximum dose and mean dose have been considered to produce clinically practical treatment plans. Objective function of the presented model includes both quadratic and linear terms, so we named the proposed model as *QuadLin*. In the next section, methods of the research have been explained; results, discussion and conclusion have formed the following sections of the paper, respectively.

## II. Methods
### A. Mathematical model

The aim in this paper is to develop a mathematical model that achieves a Pareto optimal and clinically applicable treatment plan by using dose prediction and improving it as much as possible towards the prescribed dose. In other words, the proposed model pursues two goals: 1- Achieving the predicted dose 2- Maximum improvement of the treatment plan towards the prescribed dose. In general, quadratic and linear terms in the objective function have been used to achieve the mentioned goals, respectively. The proposed objective function consists of four parts, each of which contains several terms. The four sections are: 1- Optimizing the dose of target volume voxels, 2- Controlling the dose of voxels of OARs, 3- Optimizing the maximum dose of voxels of related OARs, and 4- Optimizing the mean dose of related OARs.

The first section of the objective function as shown in equation (1) includes three terms, tow quadratic in order to minimize underdose and overdose of the PTV voxels which are denoted by non-negative variables $UD_v$ and $OD_v$, respectively, and the third one is a linear term that tries to make the delivering dose of voxels closer to the prescribed dose as much as possible. The weight of each voxel is denoted by $\omega_v$ which can also indicate the volume of the voxel, and the importance coefficient of each term in the objective function is denoted by $\psi_i$ and $\xi_j$. In general, the coefficients of terms that try to reach the predicted dose are indicated by $\psi$, and the coefficients of terms that attempt to improve the dose and reach the prescribed dose are indicated by $\xi$ which are smaller than $\psi$ in the same section.

$$min\ z_1 = \frac{\psi_1 \sum_{v \in PTV} \omega_v * UD_v^2}{\sum_{v \in PTV} \omega_v} + \frac{\psi_2 \sum_{v \in PTV} \omega_v * OD_v^2}{\sum_{v \in PTV} \omega_v} + \frac{\xi_1 \sum_{v \in PTV} \omega_v |Pres_v - \sum_{b \in B} A_{v,b} X_b|}{\sum_{v \in PTV} \omega_v} \quad (1)$$

Related constraints are shown in equations (2) and (3). $B$ is set of all beamlets $b$, each of which has an intensity of $X_b$. Dose influence matrix is denoted by $A$ and the dose delivered to voxel $v$ from beamlet $b$ with a unit intensity is shown by $A_{v,b}$. Therefore, the left side of relationships (2) and (3) is equal to the delivering dose to voxel $v$ of the PTV.

$$\sum_{b \in B} A_{v,b} X_b \geq min\{Pred_v, Pres_v\} - UD_v \qquad \forall v \in PTV \quad (2)$$

$$\sum_{b \in B} A_{v,b} X_b \leq max\{Pred_v, Pres_v\} + OD_v \qquad \forall v \in PTV \quad (3)$$

The right hand of (2) determines the amount of voxel underdose ($UD_v$) relative to the minimum of predicted dose ($Pred_v$) and the prescribed dose ($Pres_v$) of the voxel $v$. In fact, this constraint along with the first term of (1) seek to ensure that the dose delivered to the voxel is not less than the minimum of predicted and prescribed dose. The right hand of (3) determines the amount of voxel overdose ($OD_v$) relative to the maximum of predicted dose ($Pred_v$) and the prescribed dose ($Pres_v$) of the voxel $v$. Mixture of equation (3) and the second term of (1) try to ensure that the dose delivered to the voxel is not more than the maximum of predicted and prescribed dose. Therefore, the above two equations are trying to maintain the delivering dose value for the target voxels within the predicted and prescribed dose range. The third term of (1) is a linear term that tries to minimize the difference between the prescribed dose and the delivering dose of voxels of PTV as much as possible. This statement, which is a linear penalty function with a gentle slope, ensures that the dose delivered to the target voxels is as close as possible to the prescribed dose. This phrase is one of the keys to improving the treatment plan resulting from this model, compared to the predicted 3D³. Figure 1 depicts the penalty functions related to the first section of the objective function in tow situations. The figure on the left corresponds to when the predicted dose is less than the prescribed dose for a PTV voxel, and the right one corresponds to the opposite. As can be seen in the figure 1, the first section of the objective function in each case tries to achieve the range between the predicted and prescribed dose and finally meet the dose prescribed for the PTV voxels.

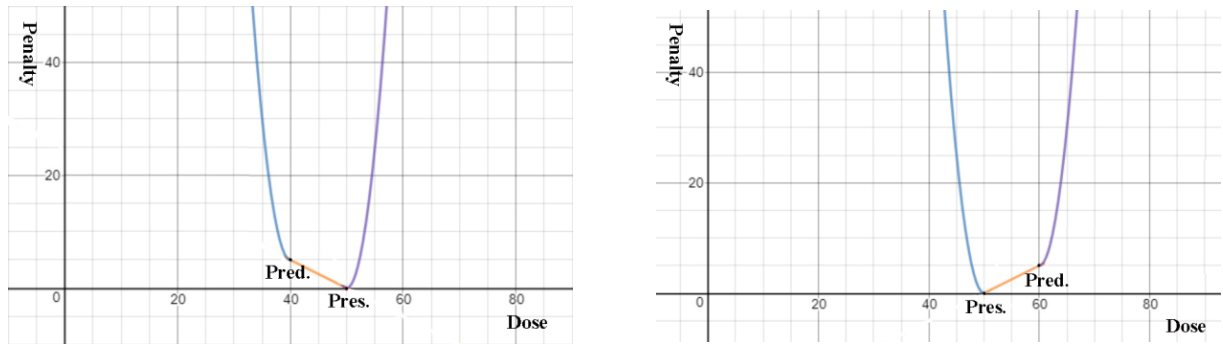

**Fig. 1. Total penalty function for voxels of the PTV (Left: predicted dos is less than prescribed dose, Right: predicted dos is more than prescribed dose)**

The second section of the objective function as shown in equation (4) includes two terms, a quadratic in order to minimize overdose of the OAR voxels which are denoted by $OD_v$, and a linear that tries to make the delivering dose of voxels closer to the zero as much as possible. The weight of each voxel as well as the importance coefficients are defined similar to the first section of the objective function.

$$min\ z_2 = \frac{\psi_3 \sum_{v \in OAR} \omega_v * OD_v^2}{\sum_{v \in OAR} \omega_v} + \frac{\xi_2 \sum_{v \in OAR}(\omega_v \sum_{b \in B} A_{v,b} X_b)}{\sum_{v \in OAR} \omega_v} \quad (4)$$

Related constraint is shown in equation (5) and notations are defined the same as section one. The left side of the relationship is equal to the delivering dose to voxel $v$ that belongs to one of the OARs.

$$\sum_{b \in B} A_{v,b} X_b \leq Pred_v + OD_v \quad\quad \forall\ v \in OAR \quad (5)$$

The right hand of (5) determines the amount of voxel overdose ($OD_v$) relative to the predicted dose ($Pred_v$) of the voxel $v$ of the OARs. This constraint seeks to ensure that the dose delivered to the voxel is not more than the predicted dose. Consequently, equation (5) tries to maintain the delivering dose value for the OAR voxels less than the predicted dose. The second term of (4) is a linear term that tries to minimize the delivering dose of voxels of OARs. Figure 2 illustrates the penalty functions related to the second section of the objective function which tries to decrease the dose of OARs voxels to predicted dose and even lower to meet the zero.

The third section of the objective function as shown in equation (6) includes two linear terms that focused on the maximum predicted dose of related OARs. In this research, limitation of the max dose of structure is implemented for three organs at risk which include brainstem, spinal cord and mandible which denoted by *OAR_Max* organs. The first statement minimizes overdose of voxels of the *OAR_Max* organs from the max predicted dose of the related structure which is denoted by $OD\_MaxP_v$, and the second term tries to decrease organs' max dose via maximizing underdose of voxels relative to related organ max predicted dose which is denoted by $UD\_MaxP_v$. The weight of each voxel ($\omega_v$) and the importance coefficient of each term in the objective function ($\psi_i$ and $\xi_j$) are defined the same as section one.

$$min\ z_3 = \frac{\psi_4 \sum_{v \in OAR\_Max} \omega_v * OD\_MaxP_v}{\sum_{v \in OAR\_Max} \omega_v} - \frac{\xi_3 \sum_{v \in OAR\_Max} \omega_v * UD\_MaxP_v}{\sum_{v \in OAR\_Max} \omega_v} \quad (6)$$

Related constraints are shown in equations (7) and (8). $OD\_MaxP_v$ and $UD\_MaxP_v$ are non-negative variables for each voxel. The left side of (7) is equal to the delivering dose to voxel $v$ of the $OAR\_Max$ organs. The right hand of (7) determines the amount of overdose of voxels of the $OAR\_Max$ ($OD\_MaxP_v$) relative to the max predicted dose of structure ($MaxP_S$). This constraint seeks to ensure that the dose delivered to the voxels is not more than the max predicted dose of related structure.

$$\sum_{b \in B} A_{v,b} X_b \leq MaxP_S + OD\_MaxP_v \qquad \forall\ v \in S,\ S \in OAR\_Max \quad (7)$$

$$\left(\sum_{b \in B} A_{v,b} X_b - \zeta_S * MaxP_S\right)^+ \leq \chi_S * MaxP_S - UD\_MaxP_v \qquad \forall\ v \in S, S \in OAR\_Max \quad (8)$$

In relationship (8) parameters $\zeta_S$ and $\chi_S$ are in the range of (0,1) which $\zeta_S$ is closer to one, $\chi_S$ is closer to zero and $\chi_S \geq 1 - \zeta_S$. In the current paper $\zeta_S = 0.9$ and $\chi_S = 0.1$, since (8) includes only non-negative values, this constraint includes the voxels of $OAR\_Max$ organs which their delivering dose is more than 0.9 of max predicted dose of related structure ($MaxP_S$), and the right side is a non-negative real number that along with the second term of (6) try to maximize the amount of voxel underdose ($UD\_MaxP_v$) relative to the 0.1 of max predicted dose of structure. In the simple word, here this constraint seeks to meet a %10 reduction in the max predicted dose of the $OAR\_Max$ structures. In general, that decline equals to $100*(1 - \zeta_S)$ percent.

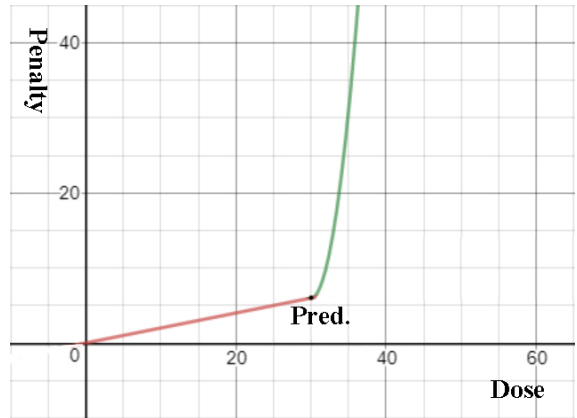

Fig. 2. Total penalty function for voxels of the OARs

The section four of the objective function as shown in equation (9) includes two quadratic terms that focused on the mean dose of related OARs. In this research, limitation of the mean dose of structure is implemented for four organs at risk which denoted by $OAR\_Mean$ organs and include right parotid, left parotid, larynx and esophagus. The first statement of (9) minimizes squared overdose of the $OAR\_Mean$ organs from the mean predicted dose of the related structure which is denoted by $OD\_MeanP_S$, and the second term tries to decrease square of organs' mean dose which is denoted by $Mean_S$. The importance coefficient of each term in the objective function ($\psi_i$ and $\xi_j$) are defined the same as the first section.

$$min\ z_4 = \psi_5 \sum_{S \in OAR\_Mean} OD\_MeanP_S^2 + \xi_4 \sum_{S \in OAR\_Mean} Mean_S^2 \quad (9)$$

Related constraints of the fourth section are as equations (10) and (11). The first one is definition of the structures' mean dose and the second constraint determines the amount of overdose of the $OAR\_Mean$ organs ($OD\_MeanP_S$) relative to the mean predicted dose of structure ($MeanP_S$). Variable $OD\_MeanP_S$ is non-negative and constraint (11) seeks to ensure that average of delivered dose to the $OAR\_Mean$ organs ($Mean_S$) is not more than the mean predicted dose of related structure ($MeanP_S$).

$$Mean_S = \frac{\sum_{v \in S} \sum_{b \in B} \omega_v (A_{v,b} X_b)}{\sum_{v \in S} \omega_v} \qquad \forall\ S \in OAR\_Mean \quad (10)$$

$$Mean_S \leq MeanP_S + OD\_MeanP_S \qquad \forall S \in OAR\_Mean \qquad (11)$$

Integrating mentioned relationships above, total mathematical model is proposed as follows:

$$\min Z = \frac{\psi_1 \sum_{v \in PTV} \omega_v * UD_v^2}{\sum_{v \in PTV} \omega_v} + \frac{\psi_2 \sum_{v \in PTV} \omega_v * OD_v^2}{\sum_{v \in PTV} \omega_v} + \frac{\xi_1 \sum_{v \in PTV} \omega_v |Pres_v - \sum_{b \in B} A_{v,b} X_b|}{\sum_{v \in PTV} \omega_v}$$

$$+ \frac{\psi_3 \sum_{v \in OAR} \omega_v * OD_v^2}{\sum_{v \in OAR} \omega_v} + \frac{\xi_2 \sum_{v \in OAR} (\omega_v \sum_{b \in B} A_{v,b} X_b)}{\sum_{v \in OAR} \omega_v} + \frac{\psi_4 \sum_{v \in OAR\_Max} \omega_v * OD\_MaxP_v}{\sum_{v \in OAR\_Max} \omega_v}$$

$$- \frac{\xi_3 \sum_{v \in OAR\_Max} \omega_v * UD\_MaxP_v}{\sum_{v \in OAR\_Max} \omega_v} + \psi_5 \sum_{S \in OAR\_Mean} OD\_MeanP_S^2 + \xi_4 \sum_{S \in OAR\_Mean} Mean_S^2$$

$S.t:$

$$\sum_{b \in B} A_{v,b} X_b \geq \min\{Pred_v, Pres_v\} - UD_v \qquad \forall v \in PTV$$

$$\sum_{b \in B} A_{v,b} X_b \leq \max\{Pred_v, Pres_v\} + OD_v \qquad \forall v \in PTV$$

$$\sum_{b \in B} A_{v,b} X_b \leq Pred_v + OD_v \qquad \forall v \in OAR$$

$$\sum_{b \in B} A_{v,b} X_b \leq MaxP_S + OD\_MaxP_v \qquad \forall v \in S, \ S \in OAR\_Max$$

$$\left( \sum_{b \in B} A_{v,b} X_b - \zeta_S * MaxP_S \right)^+ \leq \chi_S * MaxP_S - UD\_MaxP_v \qquad \forall v \in S, S \in OAR\_Max$$

$$Mean_S = \frac{\sum_{v \in S} \sum_{b \in B} \omega_v (A_{v,b} X_b)}{\sum_{v \in S} \omega_v} \qquad \forall S \in OAR\_Mean$$

$$Mean_S \leq MeanP_S + OD\_MeanP_S \qquad \forall S \in OAR\_Mean$$

$$OD_v, UD_v, OD\_MaxP_v, UD\_MaxP_v \geq 0 \qquad \forall v \in \{OAR, PTV\}$$

$$OD\_MeanP_S \geq 0 \qquad \forall S \in OAR\_Mean$$

$$X_b \geq 0 \qquad \forall b \in B$$

Finally, values of the importance coefficients of objective function's terms were allocated similar to those shown in Table 1 (M: Million, K: Kilo). Suggested model is applied on real Data and the results presented and discussed in the following sections.

**Table 1. values of the objective function's importance coefficients**

| Coefficient | $\psi_1$ | $\psi_2$ | $\psi_3$ | $\psi_4$ | $\psi_5$ | $\xi_1$ | $\xi_2$ | $\xi_3$ | $\xi_4$ |
|---|---|---|---|---|---|---|---|---|---|
| Value | 2M | 0.5M | 0.2M | 0.2M | 1K | 20K | 0.2K | 1K | 50 |

### B. Data and analysis program

In the current paper, open-access knowledge-based planning dataset (OpenKBP) [6] has been applied. The OpenKBP provides an augmented variation of real clinical data. Specifically, patient data is provided from several institutions that is available on The Cancer Imaging Archive (TCIA) [9], which hosts several open-source datasets. Current dataset stablished in order to use in the OpenKBP grand challenge (February to June 2020) and completed with the 21 series of results (i.e., dose predictions) of the challenge participants. Data for 340 patients was provided who were treated for head-and-neck cancer with 6 MV step-and- shoot IMRT in 35 fractions to satisfy the prescribed dose to the high-risk target (i.e., PTV70) from nine equispaced coplanar beams at angles 0, 40, 80, . . ., 320 degrees. The data is split into training (n=200), validation (n=40), and testing (n=100) sets. Every patient in these datasets has a dose distribution, CT images, structure masks, a feasible dose mask (i.e., mask of where dose can be non-zero), and voxel dimensions. Moreover, 21 series of dose predictions for each patient of testing sets (n=100) has been provided by the challenge participants which 19 series have better prediction quality (19 different predictions for each of the 100 patients). Therefore, complete real data for 100 patients along with 21 series of dose prediction have been shared publically in the OpenKBP dataset. Dose influence matrices were generated using the same parameters in the IMRTP library from A Computational Environment for Radiotherapy Research (CERR) [10].

In the current paper, 30 sets randomly selected from the 100 testing sets (from set 241 to 340) of the Open KBP dataset. Hence, real data of 30 patients with head and neck cancer and their 19 sets of better quality dose prediction, constitute the data of the current paper. The random patients selected via MATLAB software are sets: 245, 249, 251, 253, 254, 259, 262, 265, 270, 282, 292, 293, 294, 296, 299, 300, 301, 303, 307, 308, 309, 312, 316, 317, 319, 329, 330, 334, 335 and 339.

To solve the suggested quadratic-linear mathematical model (*QuadLin*) we used CVX, a Matlab-based convex modeling framework for specifying and solving convex programs [11], [12]. An academic CVX professional license has been applied to use CVX with commercial solver Mosek. Matlab software of version R2018a applied in this paper on a single computer with an Intel Core i5-8250U (1.80 GHz) CPU and 16 GB of RAM.

### C. Plan evaluation approaches

In order to measure the performance and quality of generated treatment plans in comparison of a reference plan, two evaluation approaches are used in this research. One of them is DVH point differences, which is the absolute difference between a DVH point from two dose distributions. The DVH differences are evaluated on two and three DVH points for each organ-at-risk (OAR) and planning-target-volume (PTV), respectively. The OAR DVH points were the $D_{mean}$ and $D_{0.1cc}$, which was the mean dose delivered to the OAR and the maximum dose delivered to 0.1cc of the OAR, respectively. The PTV DVH points were the $D_1$, $D_{95}$, and $D_{99}$, which was the dose delivered to 1% (99th percentile), 95% (5th percentile), and 99% (1st percentile) of voxels in the PTV, respectively. The second evaluation approach regarding to measure plan quality is clinical criteria satisfaction, which is defined as the proportion of clinical criteria that were satisfied by a treatment plan. This approach focuses on one criterion of each region of interest (ROI) as it can be seen in table 2.

At the end of this section, the conceptual diagram of the article is displayed in figure 3, which includes the general steps of conducting research and the relationships between them.

**Table 2. clinical criteria for each ROI**

| ROIs | Criteria |
|---|---|
| OARs | |
| Brainstem | $D_{0.1cc} \leq 50.0$ Gy |
| Spinal cord | $D_{0.1cc} \leq 45.0$ Gy |
| Right parotid | $D_{mean} \leq 26.0$ Gy |
| Left parotid | $D_{mean} \leq 26.0$ Gy |
| Esophagus | $D_{mean} \leq 45.0$ Gy |
| Larynx | $D_{mean} \leq 45.0$ Gy |
| Mandible | $D_{0.1cc} \leq 73.5$ Gy |
| PTVs | |
| PTV56 | $D_{99} \leq 53.2$ Gy |
| PTV63 | $D_{99} \leq 59.9$ Gy |
| PTV70 | $D_{99} \leq 66.5$ Gy |

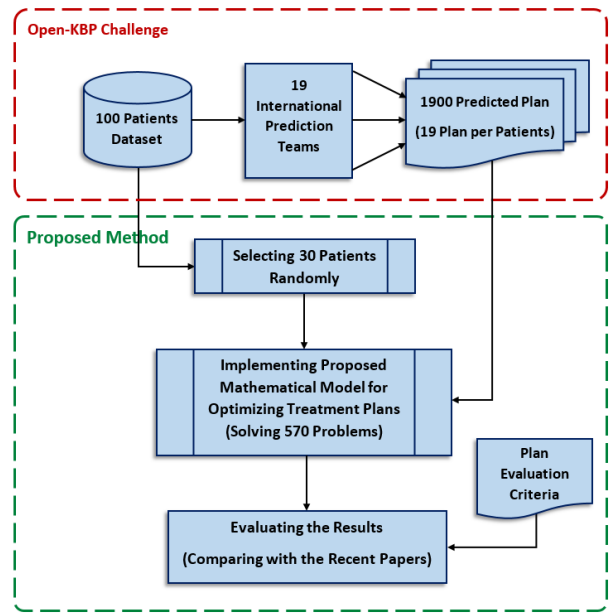

**Fig. 3. Conceptual diagram of the research**

### III.    Results

The proposed model was implemented independently on the data of 30 random patients with 19 dose predictions for each patient, consequently 570 problems were solved. Duo to evaluate the quality of plans which produced by the *QuadLin* model, two approaches were applied as defined in section II-C.

#### A.    DVH point differences

To compare the relative quality of dose distributions and review the performance of the proposed model, DVH point differences between the reference plan and both of *QuadLin*-produced dose as well as the predicted 3D$^3$ were calculated and illustrated using violin charts in figure 4. The differences were evaluated over the five DVH points listed in section II-C which lower values were better for $D_{mean}$, $D_{0.1cc}$, and $D_1$; and higher values were better for $D_{95}$ and $D_{99}$. In all five charts of the figure 4, dashed line on zero represents the reference plan; and green and red violins depict difference values for the plans of *QuadLin* and prediction, respectively. The results indicate a significant improvement in obtained treatment plans that generated by the *QuadLin* compared to both of the reference plan and 3D$^3$ predictions. The most effectiveness of the *QuadLin* can be expressed as meaningful OAR preservation (see figures 4.a and 4.b related to $D_{mean}$ and $D_{0.1cc}$ respectively) as well as desirable PTV overdose controlling (figure 4.c for $D_1$), while maintaining PTV dose and improving that as possible (figures 4.d and 4.e related to $D_{99}$ and $D_{95}$ respectively).

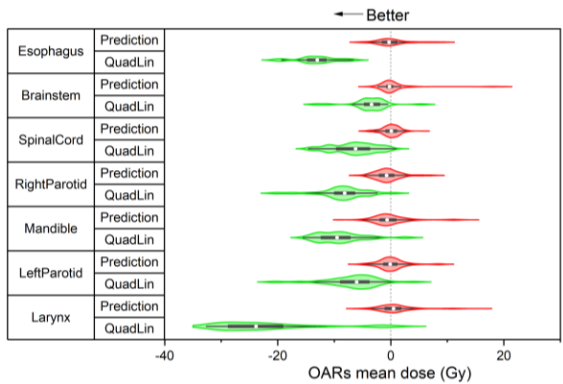
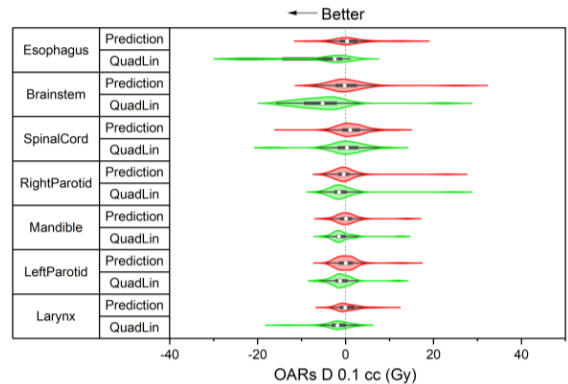

Fig. 4.a                                                                                                  Fig. 4.b

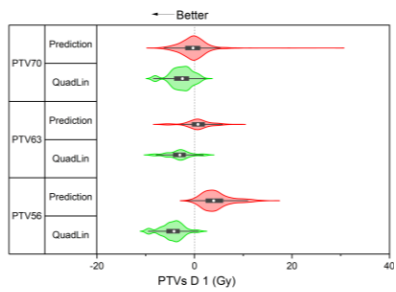
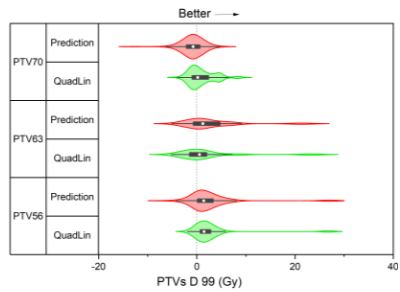
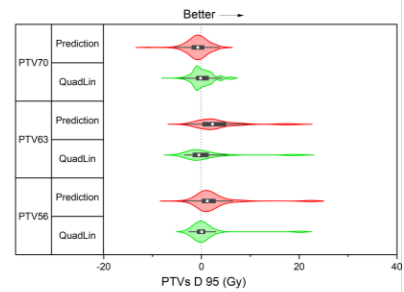

Fig. 4.c                                                    Fig. 4.d                                                    Fig. 4.e

**Fig. 4. DVH points differences between the *QuadLin*-generated plan (green violins) and 3D$^3$ prediction (red violins) with the reference plan (zero dashed lines)**

Sample DVH diagrams were drown for patient 251 applying dose prediction set 2, in figure 5, in order to compare results of the proposed *QuadLin* model (dashed lines) with the predicted 3D$^3$ (solid lines). OARs' DVH were illustrated in figures 5.a and 5.b, and targets' DVH were depicted in figure 5.c. According to the figure 5, *QuadLin* decreased the dose level of all organ at risks strongly significant, while improved both overdose and underdose of targets.

#### B.    Clinical criteria satisfaction

In this approach, proportion of the clinical criteria that were satisfied by both *QuadLin* and predicted plans for 30 patients are calculated separately for each of the 19 prediction sets. Clinical criteria consist of one criterion of each ROI as shown in table 2, which are $D_{mean}$ or $D_{0.1cc}$ for OARs and $D_{99}$ for the targets. Figure 6 displays the percentage of satisfied clinical criteria for *QuadLin* and Prediction by green and orange bullets, respectively; Horizontal dashed line represents the reference plan situation. Region of interests are categorized in OARs and targets which were illustrated by figures 6.a and 6.b, respectively; figure 6.c visualized percentage of all ROIs satisfied criteria. The results indicate that on average for 570 problems (30 patients with 19 prediction sets), presented *QuadLin* model (green bullets) improved the percentage of clinical criteria satisfaction more than 20%, 23% and 21% compared to the predicted dose (orange bullets) for OARs, targets and all ROIs, respectively.

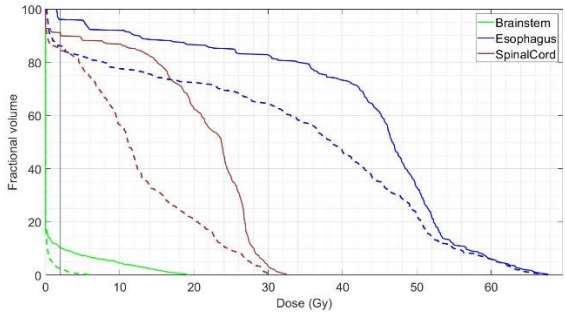
**Fig. 5.a**

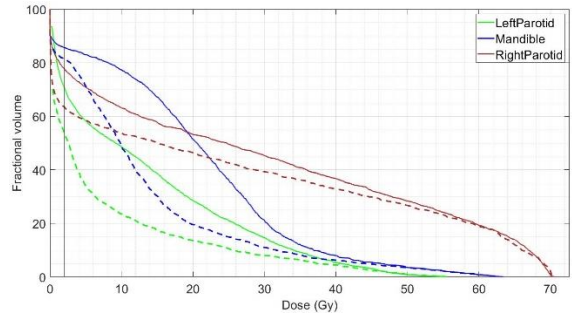
**Fig. 5.b**

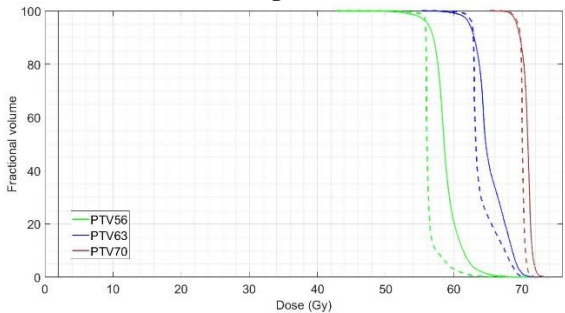
**Fig. 5.c**

**Fig. 5. Sample DVH diagram of the proposed *QuadLin* model (dashed lines) compared with the predicted 3D³ (solid lines) for patient 251 which applied prediction set 2**

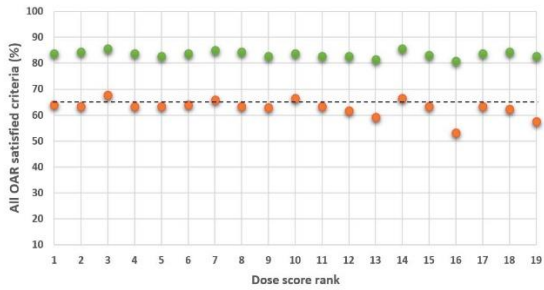
**Fig. 6.a**

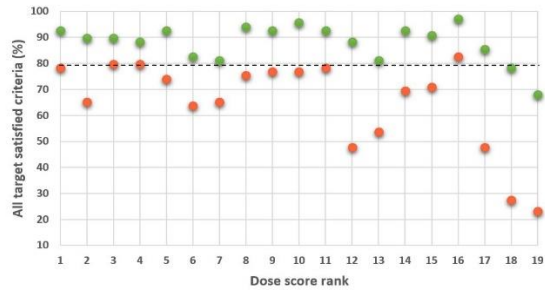
**Fig. 6.b**

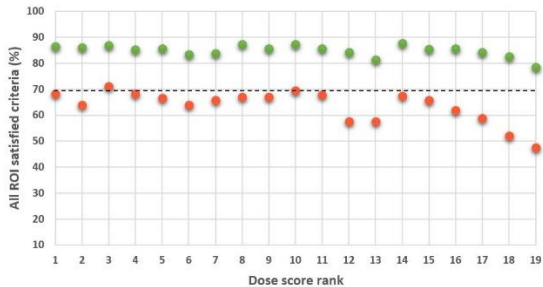
**Fig. 6.c**

**Fig. 6. Percentage of satisfied clinical criteria for OARs (fig. 6.a), targets (fig. 6.b), and all ROIs (fig. 6.c) by 3D³ prediction (orange bullets) compared with *QuadLin*-generated plans (green bullets). Dashed line depicts the reference plan situation.**

### IV.    Discussion

In recent years, extensive research has been conducted to predict three-dimensional dose distribution using artificial intelligence and deep learning tools. However, what has been less discussed is how to use and implement the predicted dose distribution, clinically and realistically. Few studies have been performed in the research literature on the extraction of an applicable treatment plan. In addition, in the same studies, the treatment plan is not well optimized and is far from the dose prescribed by oncologist. Therefore, the vacancy of a mathematical optimization model in order to extract the treatment plan and then its maximum optimization to reach the Pareto optimality level was strongly felt. *QuadLin*, the model presented in this paper, using soft constraints and integrating the predicted 3D$^3$ and the prescribed dose, as well as the objective function consisting of a combination of quadratic and linear expressions, has succeeded in presenting optimal treatment plans. Fundamental logic of the model is to first achieve the predicted dose and then improve it as much as possible.

 In this research, the valuable public dataset provided by Babier et al. [7], called Open KBP, was used, so that 30 patients were randomly sampled from 100 patients and the related data along with 19 sets of predicted dose were applied. Consequently, 570 problems were solved by the *QuadLin* model and the results were evaluated by both DVH points differences and satisfied clinical criteria approaches and compared with the reference plan as well as the predicted 3D3. The results indicate a meaningful improvement in DVH points of obtained treatment plans that generated by the *QuadLin* compared to both of the reference plan and 3D$^3$ predictions, as well as a 21% growth of the percentage of clinical criteria satisfaction compared to the predicted dose for all ROIs.

The results of the current study indicate a strong significant improvement in DVH and clinical indicators compared to previous researches. Since the model presented by Fan et al. [4] tried to make the treatment plan as close as possible to the predicted dose by using norm 2, it is clear that our model is of higher quality. Furthermore, their model only considered the predicted dose, hence, the results of their model will be seriously affected by the low quality of the prediction.

Here we compare the results of *Quadlin* model with four models presented by Babier et al. [8] titled by MaxAbs, MeanAbs, MaxRel and MeanRel. Comparison of all models is based on the results of 30 randomly selected patients (see II-B). Figure 7 illustrates average of satisfied clinical criteria for all ROIs of 30 patients resulted by each methods based on 19 dose prediction data. Horizontal black dashed line represents the reference plan. Accordingly, *QuadLin* model (green bullets) improved the clinical indicators by average more than 15%, 13%, 17% and 14% compared to the MaxAbs, MeanAbs, MaxRel and MeanRel models, respectively; and more than 21% and 15% compared to the predicted dose and the reference plan (dashed line), respectively. It is important to note that for all 19 dose prediction sets, the results of *QuadLin* model were meaningfully higher than the reference level. This is a remarkable point because it shows that reducing the quality of dose prediction can not significantly undermine the *QuadLin* results. Table 3 depicts all models' percentage of satisfied clinical criteria for each ROI. According to the table, *Quadlin* model has achieved very impressive results and improved the criteria for almost all ROIs. It satisfied approximately 85% of all clinical criteria that is by far more than all other models. The best result for each ROI is displayed with green bold.

**Table 3. Percentage of satisfied clinical criteria for each ROI resulted by each model**

|              | Reference | Prediction | MaxAbs | MeanAbs | MaxRel | MeanRel | QuadLin |
|---|---|---|---|---|---|---|---|
| *Brainstem*    | 96.6  | 97.77 | 99.39 | 99.60     | 98.99 | 99.60     | ***100.00*** |
| *Esophagus*    | 93    | 88.07 | 90.18 | ***100.00*** | 94.39 | ***100.00*** | ***100.00*** |
| *Larynx*       | 37.7  | 23.53 | 37.77 | 68.73     | 35.29 | 52.94     | ***88.24***  |
| *LeftParotid*  | 30.6  | 29.76 | 33.03 | 45.74     | 34.85 | 42.65     | ***66.73***  |
| *Mandible*     | 87.5  | 90.13 | 98.03 | 99.12     | 98.03 | 98.90     | ***100.00*** |
| *RightParotid* | 32.3  | 28.31 | 32.30 | 39.75     | 32.30 | 37.57     | ***54.91***  |
| *SpinalCord*   | 95.5  | 88.89 | 93.37 | ***99.61*** | 92.79 | ***99.61*** | 89.84    |
| *PTV56*        | 91.2  | 84.41 | ***91.50*** | 82.19 | 83.40 | 83.00     | 90.26    |
| *PTV63*        | 90.5  | 79.76 | 82.59 | 75.30     | 79.35 | 78.54     | ***86.59***  |
| *PTV70*        | 64    | 41.93 | 51.23 | 34.21     | 43.86 | 37.19     | ***86.99***  |
| *All OARs*     | 65.5  | 62.94 | 67.95 | 76.68     | 68.23 | 74.13     | ***83.55***  |
| *All Targets*  | 79.4  | 65.06 | 72.31 | 60.03     | 65.45 | 62.24     | ***88.16***  |
| *All ROIs*     | 69.7  | 63.56 | 69.22 | 71.81     | 67.42 | 70.65     | ***84.89***  |

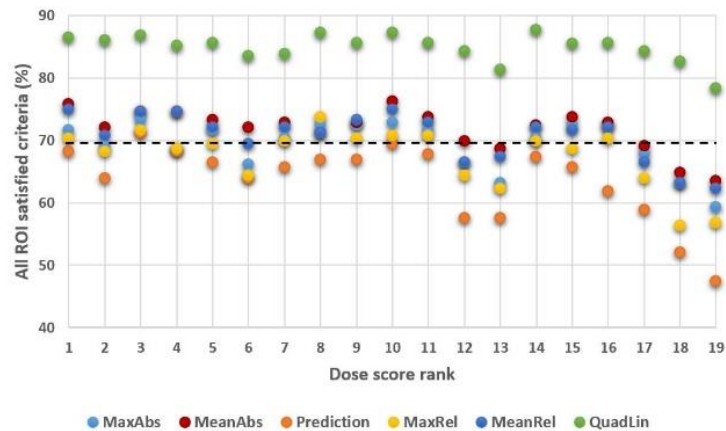

**Fig. 7. Satisfied clinical criteria for all ROIs of 30 patients in each 19 dose prediction sets; comparing all models**

### V.     Conclusion

In the current research, an attempt was made to revolutionize the knowledge-based planning by presenting a new mathematical model, and to take a serious step towards optimizing the treatment plan derived from predicted three-dimensional dose distribution. The CVX framework with commercial solver Mosek as well as valuable open access dataset, Open-KBP, helped us a lot in this research. Evaluation of the research results indicates the significant effectiveness of the *QuadLin* model in that while reducing the dose and preserving organs at risk, it has also improved the dose delivery to the target volumes. Based on the literature, this article has generated the best treatment plan from the predicted 3D[3] so far.

**Acknowledgments**
The first author would like to express his respect, appreciation and gratitude to his old friend and master who played a significant role in shaping this article; a friend who wished to remain anonymous.

**Conflict of Interest Statement**
The authors have no relevant conflicts of interest to disclose.